\documentclass[12pt]{article}
\usepackage{amsmath,amsfonts,amssymb,amsthm}

\def\noi{\noindent}

\def\pf{\noi{\bf Proof.\ \,}}
\def\eop{{$\square$}}

\def\labtt#1{\label {#1} }

\def\refpp#1{(\ref {#1})}

\def\a{\alpha}
\def\b{\beta}
\def\g{\gamma}

\def\l{\lambda}

\def\s{\sigma}

\def\CC{{\mathbb C}}
\def\FF{{\mathbb F}}

\def\RR{{\mathbb R}}

\def\ZZ{{\mathbb Z}}
\def\MM{{\mathbb M}}

\def\tr{{\rm tr}}

\def\la{\langle}
\def\ra{\rangle}
\def\<{\langle}
\def\>{\rangle}

\def\bs{\it}            

\def\dim{{\bs dim}}

\def\exp{{\bs exp}}

\def\l{{\lambda}}

\def\dual#1{#1^*}        

\def\vac{\hbox{\bf 1}} 

\def\vnat{{V^\natural}}

\begin{document}

\newtheorem{thm}{Theorem}[section]
\newtheorem{prop}[thm]{Proposition}
\newtheorem{lem}[thm]{Lemma}
\newtheorem{rem}[thm]{Remark}
\newtheorem{coro}[thm]{Corollary}
\newtheorem{conj}[thm]{Conjecture}
\newtheorem{de}[thm]{Definition}
\newtheorem{hyp}[thm]{Hypothesis}

\newtheorem{nota}[thm]{Notation}
\newtheorem{ex}[thm]{Example}
\newtheorem{proc}[thm]{Procedure}

\newcommand{\VAg}{V\!\!A((V_X^\gamma)_1)}

\newcommand{\IVAg}{V\!\!A((IV_X^\gamma)_1)}

\newcommand{\IVg}{IV(X)^\gamma}

\begin{center}

{\Large \bf
Applications of vertex algebra covering procedures to Chevalley groups and modular moonshine}

\vspace{10mm}
Robert L.~Griess Jr.
\\[0pt]
Department of Mathematics\\[0pt] University of Michigan\\[0pt]
Ann Arbor, MI 48109  USA  \\[0pt]
{\tt rlg@umich.edu}\\[0pt]
\vskip 1cm

Ching Hung Lam
\\[0pt]
Institute of Mathematics \\[0pt]
Academia Sinica\\[0pt]
Taipei 10617, Taiwan\\[0pt]
{\tt chlam@math.sinica.edu.tw}\\[0pt]
\vskip 1cm
\end{center}

\begin{abstract}
A vertex operator algebra of lattice type ADE has a standard integral form  which extends a Chevalley basis for its degree 1 Lie algebra.   This integral form may be used to define a vertex algebra over a commutative ring $R$ and to get a Chevalley group over $R$ of the same type, acting as automorphisms of this vertex  algebra.     We define vertex algebras of types BCFG over a commutative ring and certain reduced VAs, then get analogous results about automorphism groups.    In characteristics 2 and 3, there are exceptionally large automorphism groups.   A covering algebra idea of Frohardt and Griess for  Lie algebras is applied to the vertex algebra situation.

We use integral form and covering procedures for vertex algebras to complete the modular moonshine program of Borcherds and Ryba for proving an embedding of the sporadic group $F_3$ in $E_8(3)$.
\end{abstract}

\newpage

\tableofcontents

\newpage

\begin{center}
{\bf \Large Table of Notations}
\vskip 0.1cm
\small
\begin{tabular}{|c|c|c|}
   \hline
\bf{Notation}& \bf{Explanation} & \bf{Examples }  \cr
&& {\bf in text}\cr
   \hline\hline
ancestor of $X'$  & $(X,\g)$ for certain root system& \refpp{exceptionalpair}\cr
& and graph automorphism& \cr
\hline
$A_1,  \cdots , E_8$ & root system or root lattice  &  \refpp{defvanonsimplylaced} \cr
& of type $A_1, \dots , E_8$& \cr
\hline
classical Lie algebra & has a Chevalley basis & \refpp{classicalliealg}\cr
\hline
classical VA & $R \otimes IV_X = R\otimes V(X)$ ($X$ of type $A, D, E$)   & \refpp{v(y)}, \refpp{classicalva}\cr
vertex algebra&
 or $R \otimes IV(Y)$ ($Y$ of type $B, C, G, F$)                 & \refpp{standardif}  \cr
\hline
covering algebra, CA & An subalgebra (or subVA) $B$ of   & \refpp{coveringsubgroup}, \cr
covering VA  & an algebra (or VA) $A$ such that &  \cr & $A=I +  B$ for an ideal $I$&  \cr
 \hline
 covering subgroup & a subgroup $C$ that maps onto $A/B$,  & \refpp{coveringsubgroup}\cr
of $A$ w.r.t $B$ & i.e., $A=B+C$ & \cr
\hline
covering transversal &   A covering subgroup $C$     & \refpp{coveringsubgroup}\cr
& such that $A=B\oplus C$ &\cr
\hline
$EE_8$ & lattice isometric to $\sqrt 2$ times the $E_8$ lattice&  \refpp{not6.1}  \cr
\hline
exceptional pair & certain $(X',p)$: (root system, prime) & \refpp{exceptionalpair}, \refpp{fg1} \cr
\hline
$G(X, R)$ & the Chevalley group on the Lie algebra  &  \refpp{chevgrouponva1}\cr
& with Chevalley basis of type $X$ over $R$& \refpp{chevgroup}\cr
\hline
$G(X, R)^\gamma$ &  the fixed point subgroup of $\gamma$ in $G(X,R)$ & \refpp{chevgroup} \cr \hline
IF, subIF & integral form, sub integral form & \refpp{chevgroup}\cr
\hline
$\nu$ & the norm element of a finite group, &  \refpp{normideal}, \cr
      & i.e., $\nu= \sum_{g\in G} g \in \ZZ[G]$ &  \refpp{H0} \cr
\hline
norm ideal & the image of the norm map &  \refpp{normideal}\cr
\hline
reduced algebra & the quotient of the fixed point subalgebra &  \refpp{normideal} \cr
& by the norm ideal& \cr\hline
$s_{\alpha, n}$ & the coefficient of $z^n$ in $\exp( \sum_{n>0} \frac{\alpha(-n)}{n} z^n)$&  \refpp{covering}\cr
& i.e., $\exp\left( \sum_{n>0} \frac{\alpha(-n)}{n} z^n\right) = \sum_{n\geq 0} s_{\alpha,n}z^n$ & \cr \hline
VA, VOA & vertex algebra, vertex operator algebra & \refpp{defva},\refpp{defvoa} \cr
\hline
$V\!\!A(W)$ & the subVA generated by a subset $W$ &  \refpp{defvanonsimplylaced} \cr \hline
$V_L$ & lattice VOA associated to $L$ &  Section \ref{sec:2} \cr\hline
$IV_L$ or $IV(L)$ & the standard integral form of $V_L$ &  \refpp{standardif}\cr
\hline
$V_L[1/2]$ & the $\ZZ[1/2]$-form  $IV_L\otimes \ZZ[1/2]$ of   $V_L$ & \refpp{12form}\cr \hline
$V^\natural[1/2]$ & a $\ZZ[1/2]$-form  of   $V^\natural$ & \refpp{p6.7} \cr \hline
${^gV}$ & the vertex algebra $ (V^\natural[1/2])^g/ \nu(V^\natural[1/2])$ &  \refpp{defgv}\cr
& associated to an odd order element $g\in \MM$ & \cr \hline
$Z( \cdot )$ & center of an algebra or group & \refpp{fglem2}\cr \hline
 \end{tabular}
\end{center}
\newpage

\section{Introduction}

Recent results on integral forms in vertex algebras \cite{borcherds}, \cite{ivoa}
suggested a look at Chevalley groups acting on vertex algebras over arbitrary commutative rings.   We begin this article by explaining how to define classical vertex algebras for all types of root systems, and Chevalley groups acting on them.

As with modular Lie algebras (i.e., in positive characteristic), there is exceptional behavior for certain types of classical vertex algebras in characteristics 2 and
3.    Such a  reduced vertex algebra (the vertex algebra modulo a certain nontrivial ideal) \refpp{normideal}  has automorphism group which is  a larger Chevalley group.

\begin{de} \labtt{coveringsubgroup}
Given an abelian group $A$ and subgroup $B$, the subgroup $C$ is called a covering subgroup with respect to $B$ if $C$ maps onto $A/B$, i.e.,
$A=C+B$.     A covering transversal is a covering subgroup $C$ so that $A=B\oplus C$.

If $A$ is an algebra, respectively, a vertex algebra,  and $B$ is an ideal,  a subalgebra, respectively, a subVA, $C$,  is said to cover $A/B$ if $C+B=A$. $C$ is
called a covering algebra or a covering VA with respect to $B$.
\end{de}

The covering procedure developed by Frohardt and Griess
\cite{fg}, based on action of a graph automorphism, was used to demonstrate this exceptional behavior for Lie algebras in a uniform  way (linear algebra of a graph automorphism) without case-by-case work and calculations used in earlier treatments (see \cite{fg} for details and history).

Fortunately, these covering procedures can easily be promoted to the vertex algebra situation to construct  the exceptional automorphism groups for corresponding vertex algebras.

\begin{ex} \labtt{a2g2}
Here is an example for Lie algebras in characteristic 3.   The algebra $a_2(F)$ is 8-dimensional.   The quotient  $a_2(F)$ modulo its 1-dimensional center, is simple.   Because of  embeddings $a_2(F) < g_2(F) < d_4(F)$ defined by a graph automorphism of order 3, we can deduce  that the automorphism group of the reduced $a_2(F)$ is $G_2(F)$.   The subalgebra $a_2(F)$ in $g_2(F)$ is referred to as a {\it covering algebra} in \cite{fg}.
For the analogue in vertex algebra theory, we consider a certain containment of integral forms $ IV_{A_2} <  IV(G_2) <  IV_{D_4}$ (see \refpp{XSg}) defined by a graph automorphism,  then argue that a quotient of $F \otimes IV_{A_2}$ by an ideal has automorphism group $G_2(F)$.
\end{ex}

A variant of the  covering procedure \cite{fg} applies to the interesting study by Borcherds and Ryba of a $3C$ element $g$  in the Monster simple group.   Its centralizer has the form $C(g) \cong 3 \times F_3$, where $F_3$ is a sporadic simple group of order $2^{15}3^{10}5^3 7^2 13{\cdot }19{\cdot} 31$.   Computer work with a 248-dimensional lattice proved that $F_3$ embeds in the Chevalley group $E_8(3)$ \cite{th}.   Borcherds and Ryba created a graded $\FF_3$ vector space which was a module for $C(g)/\la g \ra$ and which they felt ought to be (up to some re-indexing)
$IV_{E_8}/3 IV_{E_8}$.   We verified their conjecture by use of several vertex algebras and a subspace of one which plays the role of a covering transversal \refpp{coveringsubgroup}.

Using the covering procedure, one might devise new embedding proofs of finite groups into finite groups of Lie type and other finite groups.
Besides the standard integral forms $IV_L$, one might consider examples in \cite{ivoa}, which include a Monster-invariant integral form in the Moonshine VOA.

\section{VAs and VOAs over commutative rings}

We begin by reviewing some basic concepts.

\begin{de}\labtt{defva}
Let $R$ be a commutative ring with identity.  A vertex algebra (VA) over $R$  is a $R$-module $V$ equipped with a linear map
\begin{eqnarray*}
Y(\ ,z):V &\longrightarrow &\mathrm{End}\ V\left[ \left[ z,z^{-1}\right]
\right] \\
v &\longrightarrow &Y(v,z)=\sum_{i\in \mathbb{Z}}v_{i}z^{-i-1}
\end{eqnarray*}
and a linear map $T:V \to V$ such that the following conditions hold:

\begin{enumerate}
\item  (Vacuum condition) there is a vector $\vac$ such that
\[
Y(\vac,z)=id_V = id_V \cdot z^0;
\]

\item  (Creation property) $Y(a,z)\cdot \vac\in V\left[ \left[ z\right] \right]
$ and $a_{-1}\cdot \vac =a$ for any $a\in V$;

\item  for any $a,b\in V$,  $a_{n}b=0$  for  $n$  sufficiently large;

\item  (Translation  property)   $T\vac =0$ and
\[
[T, Y( a,z)] =Y(Ta, z)= \frac{d}{dz}Y(a,z);
\]

\item  (Borcherds Identity) for any $a,b,c\in V$ and  $m,n, q\in \mathbb{Z}$,
\[
\begin{split}
&\sum_{i\in \mathbb{Z}} \binom{m}{i} (a_{q+i}b)_{m+n-i} c\\
 = & \sum_{i\in \mathbb{Z}} (-1)^ i  \binom{q}{ i} \left( a_{m+q-i} (b_{n+i}c) -  (-1)^q  b_{n+q-i} (a_{m+i} c) \right).
 \end{split}
 \]
\end{enumerate}

A vertex algebra $V$ is $\mathbb{Z}$-graded if  $ V=\oplus _{n\in \ZZ}  V_{n}$ such that  if  $a\in V_n$, $b\in V_m$ are
homogeneous, then $a_k b\in V_{n+m-k-1}$ for any integer $k$.
\end{de}

\begin{de}\labtt{defvoa}
A $\ZZ$-graded VA  $V=\oplus _{n\in \ZZ}  V_{n}$ over $R$ is said to be a  \emph{vertex operator algebra (VOA)}  if there is a special
element $\omega \in V$ such that the operators $L(n)=\omega_{n+1}, n\in \ZZ,$ satisfy the Virasoro relation
\[
[L(m) , L (n)] = (m - n) L (m + n) + \binom{m+1}3 \delta_{m+n,0} c'
\]
for  any $m, n\in \ZZ$, where $c' \in R$.   We call  $2c'$ the rank or the central charge of $V$.  (This definition of rank is worded  differently from the usual one since $R$ may not contain $\frac 12$.)    Moreover,
\begin{enumerate}
\item  the $\ZZ$-grading of $V$ is compatible with the action of $L(0)$, i.e.,
\[
L (0) v = nv = wt (v) \cdot v\quad \text{ for any }v \in V_n;
\]
(note:  $n\in \ZZ$, whereas $L(0) \in \ZZ\cdot 1_R$

\item $\dim V_n< \infty $ for any $n\in \ZZ$ and $V_n=0$ for $n$ sufficiently negative;

\item $T=L(-1)$. In particular,  $\frac{d} {dx} Y (v, z) = Y (L ( -1) v, z)$ for any $v\in V$.
\end{enumerate}
\end{de}

\begin{de}\labtt{defif} Suppose that $V$ is a VOA (over the complex numbers) with a
nondegenerate symmetric bilinear form. An integral VOA form (abbreviated
IVOA) for $V$ is an abelian subgroup $J$ of $(V, +)$ such that $J$ is a VA over
$\mathbb{Z}$, there exists a positive integer $s$ so that $s\omega\in J$, for each $n$, $J_n := J\cap V_n$
is an integral form of $V_n$, $(J, J) \in \mathbb{Q}$. Since $J$ is a VA, $\mathbf{1}\in J$. For each degree $n$, $J_n$ has finite rank, whence there is an integer $d(n) > 0$ so that
$d(n)\cdot (J_n, J_n)\in \mathbb{Z}$.
\end{de}

\begin{rem}
An integral form $J$ of a VOA over $\CC$ will be a VA over $\mathbb{Z}$. If $R$ is a commutative ring, then $J\otimes_\ZZ
R $ is a VA over $R$.
\end{rem}

\section{The Chevalley basis and Chevalley group construction for vertex algebras of types ADE}\label{sec:2}

For an even integral lattice, let $V_L$ be the usual lattice vertex operator algebra.   Chevalley group notation is taken from \cite{carter}.

\begin{nota}\labtt{aboutx}  For convenience, we may use a common symbol, $X$, to indicate a type ADE root lattice, the root system or the name of the root system.   For the case of root system of type B, C, F, G, we shall use a symbol $Y$ for the root system or its name, but it shall not refer to a lattice.
\end{nota}

One reason for the restriction on $Y$ is that lattice type VOAs require the defining lattice to be even.

\begin{nota}\labtt{standardif}   When $V=V_L$, denote the standard integral form by  $IV_L$ or $IV(L)$.   It is given in \cite{borcherds}, \cite{ivoa}.   (A variation of this is positive definite.)
\end{nota}

Let  $L$ be a root lattice.  This standard $IV_L$ admits the group of graph automorphism of the root system,  lifted in obvious way from $L$ to $V_L$ (see \cite{gh}, Appendix).

If $(IV_L)_n:=IV_L\cap V_n$, then $(IV_L)_1$ is a lattice in the Lie algebra $(V_1, 0^{th})$ which is  spanned by a Chevalley basis, $\mathcal C$.
In \cite{carter},  notation for elements of $\mathcal C$ in root spaces for root $\a$ is $e_\a$.
If $t \in \CC$, then $x_\a (t):=\exp(t\, ad(e_\a))$ takes
$\ZZ [t] \otimes (IV_L)_1$ to itself.   See \cite{carter}.

Since $(IV_L)_1$ generates $IV_L$ \cite{ivoa}, $x_{\a} (t)$ takes $\ZZ [t] \otimes IV_L$ to itself.
The result is that for $t\in A$,  $x_{\a} (t)$ gives an automorphism of the VOA $A\otimes IV_L$.

\begin{thm}\labtt{chevgrouponva1}
Let $X$ be  an indecomposable root system of type ADE and let $R$ be a commutative ring.   Then the natural action of the Chevalley group $G(X, R)$ on the Lie algebra with Chevalley basis of type $X$ over the scalars $R$ extends
uniquely to an action as automorphisms on the vertex algebra $R\otimes IV_L$.
\end{thm}

\section{Definition of VA of type BCFG}

\begin{de}\labtt{H0}
Let $S$ be a subset of a group acting on a group $A$. Define $A^S$ to be the elements of $A$ fixed by $S$.

In case $G:= \la S \ra$ is a finite group and $(A, +)$ is a $G$-module, we write  $\nu$ for
the norm map defined by $\nu (a ) = \sum_{g\in G} g(a )$.
The quotient $A^G/\nu (A)$ is often called the 0-th Tate cohomology group, $\hat{H}^0(G, A)$.
\end{de}

As explained in the last section, a VOA of type ADE (one root length) is just $V_L$, where $L$ is the root lattice.   For indecomposable root systems with more than one length, more care is needed.   We use the definition below for a VOA of type $B_n, C_n, G_2, F_4$.

\begin{de}   \labtt{defvanonsimplylaced}
Let $X$ be the root lattice of type ADE  and let $\g$ be the standard lift of the graph automorphism of $X$ whose fixed points
on the Lie algebra at degree 1 is a Lie algebra of type $Y$. The subVA $\VAg < V_X^\gamma$ generated by $(V_X^\gamma)_1$ is
defined to be the classical VA of type $Y$.   The cases are listed here:
\[
\begin{array}{|c|c|c|c|}
\hline
X &  Y &\text{central charge of } & \text{central charge of }\cr
& & V_X^\gamma \text{ over } \CC\ & \VAg \text{ over } \CC\cr \hline
 D_{n+1} & B_n  & n+1 &n+  \frac{1}2\cr
 A_{2n-1}&   C_n &  2n-1 &  2n -\frac{3n}{n+2}\cr
D_4 & G_2 &  4 &  \frac{14}5\cr
E_6 &  F_4  &  6 &    \frac{26}5\cr
\hline
\end{array}
\]
\end{de}

\begin{rem}
The central charge of $\VAg$ over $\CC$ is computed by the formula
$ c=k\, \dim (\mathfrak{g}) /(k+h^\vee),$
where $k$ is the level,
$\mathfrak{g}=(\VAg)_1$ and $h^\vee$ is the dual Coxeter number\cite{FZ,Kac}.

Note also the distinction between
$\VAg$ and $ V_X^\gamma$.
\end{rem}

\begin{nota}\labtt{XSg}
$X$ ADE root lattice, $\g$ graph automorphism
so that the fixed points on the Lie algebra has type $Y$;
$X^\gamma$, the fixed point sublattice of $\gamma$ in $X$. Then, $X^\gamma$ is a root sublattice
of $X$.
\[
\begin{array}{|c|c|c|c|}
\hline
X &  Y &   X^\gamma & \text{order of } \gamma \cr \hline
 D_{n+1} & B_n  & D_n &   2 \cr
 A_{2n-1}&   C_n &  A_1\oplus \cdots \oplus A_1  &  2\cr
D_4 & G_2 &  A_2 &  3\cr
E_6 &  F_4  &  D_4 &   2\cr
\hline
\end{array}
\]

In the integral form $IV_X$, we have two subobjects:

$IV_X^{\g}$, the fixed points;

$\IVAg$, the vertex subalgebra of $IV_X^{\g}$ generated by the degree 1 component $(IV_X)^{\g}_1$ of the standard integral form of
$(V_X)^{\g}_1$
 the Lie algebra at degree 1.
Note that $\IVAg$ is an integral form of $\VAg$ and
we call $\IVAg$ the standard integral form of $\VAg$.
\end{nota}

\begin{nota}\labtt{v(y)}
We use $V(Y)$ to denote a  VOA of type $Y$ over the complex numbers \refpp{defvanonsimplylaced}, i.e., $V(Y)$ is the
lattice VOA $V_Y$ if  $Y$ is of type ADE and $V(Y)= \VAg$ if $Y$ is of type BCFG and the degree 1 space $(V_X^\g)_1$ is a Lie
algebra of type $Y$.   The standard integral form in $V(Y)$ is denoted $IV(Y)$.
\end{nota}

\begin{de}\labtt{classicalva}
We define a classical vertex algebra over the commutative ring $R$ to be
$R \otimes IV_X=R\otimes IV(X)$ (for $X$ of type $A, D, E$) and
$R \otimes IV(Y)$ (for $Y$ of type $B, C, G, F$).
See
\refpp{defvanonsimplylaced}, \refpp{XSg}.
\end{de}

Thus,  a classical vertex algebra over $R$ has a ``Chevalley basis''.

\subsection{The Chevalley basis and Chevalley group construction for vertex algebras of types BCFG}

\begin{nota}\labtt{chevgroup}  The Chevalley group of type $X$ over the commutative ring $R$ is denoted $G(X,R)$.
It is, by definition, a subgroup of the Lie algebra of type $X$ over $R$ obtained from the Chevalley basis.

If $X$ is a root lattice of type ADE, the group associated to  $R \otimes IV_X$ is denoted $G(IV_X,R)$. The restriction of
$G(IV_X,R)$ to $R\otimes (IV_X)_1$ gives an isomorphism of groups $G(IV_X, R)  \rightarrow G(X,R)$.    We shall usually identify
these groups.

If $\g$ is a graph automorphism lifted in the obvious way to the group $G(X, R)$, then the fixed point subgroup $G(X,R)^{\g}$ leaves invariant the fixed point subIF $(R\otimes IV_X)^{\g}$, whence also $R\otimes
\IVAg$.
\end{nota}

\begin{thm}\labtt{chevgrouponva2}
Let $X$ be a root  system of type ADE and let $R$ be a commutative ring. Let $\g$
be a graph automorphism. Then the natural action of the Chevalley group $G(X, R)^{\g}$ on the Lie algebra with Chevalley basis of
type $Y$ (i.e., $B, C, F, G$) over the scalars $R$ extends to an action as automorphisms on the vertex algebra $R\otimes \IVAg$
\refpp{defvanonsimplylaced}.
\end{thm}

To summarize: a Chevalley group of any type over a  commutative ring $R$ acts as automorphisms of a VA which  extends the action on the Lie algebra over $R$ with Chevalley basis.

\section{The covering algebra procedure lifted to VAs}

The {\it covering algebra procedure } or {\it CA-procedure} for describing the exceptionally large automorphism groups was introduced in \cite{fg}.   See Appendix A for a summary.
It gave a new and relatively computation free determination of the exceptionally large automorphism groups.

In this section we show that the CA procedure can be promoted to our studies of VAs and automorphism groups.

\bigskip

\begin{de} \labtt{normideal}
Suppose $\g$ in \refpp{XSg} had order $p$.   Define the norm element $\nu := \sum_{j=0}^{p-1} \gamma^j \in \ZZ \la \g \ra$. Let
$Q$ be any $\g$-invariant subring of $V(X)$ (such as $V(X)$ or $IV(X)$).  Define $N:=\nu (Q)$. Then $N$ is an ideal in $Q^{\g}$,
the fixed point subring \refpp{fglem1}. We call $N$ the norm ideal of $Q$. The reduced ring, or reduced algebra, is the quotient
$Q^\g/N$ of the fixed points $Q^{\g}$ by the norm ideal.   (This is an example of a 0th Tate cohomology group \refpp{H0}.)

If $R$ is a commutative ring, the above definitions extend to $R\otimes IV(X)$.

If $(X',p)$ is exceptional \refpp{exceptionalpair} and $F$ is a field of characteristic $p$, the Lie algebra or VA of type $X'$ over $F$
has a norm ideal, defined by use of the ancestor $(X, \g)$ \refpp{exceptionalpair} of $(X',p)$ and by use of the norm element on
the Lie algebra or VA of type $X$.
\end{de}

\begin{lem}\label{covering}   Let $X, \g$ be as in \refpp{chevgroup}.
Let $IV_X$ and $IV_{X^\g}$ be the standard integral forms for the respective lattice VOA $V_X$ and $V_{X^\gamma}$. Let $N:=\nu(IV_X)$ be
the norm ideal defined in \refpp{normideal}. Then,
\[
(IV_X)^\gamma = IV_{X^\g} + N.
\]
That is, $IV_{X^\g}$ is a covering subIF  of $(IV_X)^\gamma$ with respect to $N$ \refpp{coveringsubgroup}.
\end{lem}

\pf
For any $\alpha \in X$, let $s_{\alpha, n}$ be the coefficient of $z^n$ in
$E^-(-\alpha, z) = \exp( \sum_{n>0} \frac{\alpha(-n)}{n} z^n)$, i.e,
\[
\exp\left( \sum_{n>0} \frac{\alpha(-n)}{n} z^n\right) = \sum_{n\geq 0} s_{\alpha,n}z^n.
\]
By \cite[Theorem 3.3]{ivoa},  the standard integral form $IV_X$ is spanned over $\ZZ$  by
\[
\{ s_{\alpha_1,n_1}s_{\alpha_2,n_2} \cdots s_{\alpha_k,n_k}\otimes e^\alpha \mid \alpha_i, \alpha\in X , n_i >0\}.
\]
Therefore,
\[
\begin{split}
(IV_X)^\gamma  = &\ span_\ZZ  \{ s_{\alpha_1,n_1}s_{\alpha_2,n_2} \cdots s_{\alpha_k,n_k}\otimes e^\alpha \mid \text{ all } \alpha_i, \alpha\in X^\gamma, n_i >0 \} \\
                & + span_\ZZ \{ \sum_{j=0}^{p-1} \gamma^j(s_{\alpha_1,n_1}s_{\alpha_2,n_2} \cdots s_{\alpha_k,n_k}\otimes e^\alpha)  \mid \text{ some } \alpha_i \text{ or }\alpha\notin X^\gamma \}.\\
=& IV_{X^\g} + \nu(IV_X)
\end{split}
\]
as desired.
\eop

\begin{rem}\label{rem4.4}
Let $R$ be a commutative ring. Then the automorphism $\gamma$ also acts on $R\otimes IV_X$. By the same proof as in \refpp{covering}, we have
\[
(R\otimes IV_X )^\gamma = R\otimes IV_{X^\gamma} + \nu(R\otimes IV_X).
\]
Note also that $N= \nu(R\otimes IV_X) = (R\otimes IV_X)^\gamma$ if $(char(R),|\gamma|)=1$.
\end{rem}

\medskip

Since $X^\g$ is also a root lattice, the integral form $IV_{X^\g}$ is generated by $(IV_{X^{\g}})_1$ (see \cite{ivoa}), which is
contained in $(IV_X^\g)_1$. Hence, $IV_{X^\g} < \IVAg$. By \refpp{covering} and \refpp{rem4.4}, we also have the
following lemma.

\begin{lem}\label{lem:4.5}
Let $R$ be a commutative ring. Then
\[
(R\otimes IV_X )^\gamma = R\otimes \IVAg + \nu(R\otimes IV_X).
\]
In particular, $R\otimes \IVAg$ is a covering algebra \refpp{coveringsubgroup} in $(R\otimes IV_X )^\gamma$ with respect to the
norm ideal $\nu(R\otimes IV_X)$.
\end{lem}

\subsection{More quotient algebras}
The group $G(X,R)^{\g} $ also acts as automorphisms of the quotient algebra $(R\otimes IV_X)^{\g} /  (R\otimes IV_X)^{\g} \cap
\phi(\g) ( R \otimes IV_X   ) $, for any polynomial $\phi (t) \in R[t]$ \refpp{fglem2}.

For the case $\phi (t)=(t-1)^{p-1}$, and $R$ is a field of characteristic $p:=|\g|$, then we have $\phi (\g ) =\sum_{j=0}^{p-1}
\gamma^j $ and $ \phi(\g)( R\otimes IV_X) =\nu(R\otimes IV_X) < (R\otimes IV_X)^{\g}$. By \refpp{covering}, \refpp{rem4.4} and \refpp{lem:4.5}, we have the following proposition.

\begin{prop}\labtt{largeaut}
Let $X, X^\g, \g$ and $\IVAg$ be defined as in Notation \refpp{XSg}. Let $R$ be a commutative ring. Then
\[
\frac{R\otimes IV_{X^\g} }{\nu(A\otimes IV_X)\cap (R\otimes IV_{X^\g})} \cong \frac{(R\otimes IV_X)^{\g}}{ \nu(R\otimes IV_X)}
\]
and
\[
\frac{R\otimes \IVAg}{\nu(R\otimes IV_X)\cap R\otimes \IVAg} \cong \frac{(R\otimes IV_X)^{\g}}{ \nu(R\otimes IV_X)}
\]
as VAs. In particular, the group $G(X,A)^{\g} $ acts as automorphisms of the quotient algebra $ R\otimes IV_{X^\g}/ \nu(R\otimes
IV_X)\cap (R\otimes IV_{X^\g})$.
\end{prop}

\section{Automorphism group of a classical VA over a field}

\begin{prop}\labtt{autclassicalliealg}
Let $\frak g$ be a Lie algebra over the field F.   Assume that $\frak g$ has classical type.  Then one of the two following applies.  (Here, an automorphism of a Lie algebra is required to be $F$-linear.)

\begin{enumerate}
\item[(1)] Suppose that $\frak g$ is not exceptional.  If  $Y$ is the type of  $\frak g$, then
$Aut(\frak g)$ is isomorphic to the subgroup of $Aut(G(Y,F))$ generated by inner automorphisms, and all diagonal and graph automorphisms.

\item[(2)]  Suppose that $\frak g$
is exceptional.  Then
(in the notation of
\refpp{XSg}) the reduced Lie algebra $\overline {\frak g}$ and $Aut(\overline {\frak g})$ occur in one of conclusions listed in
 \refpp{fg1}.
  We list only the cases of \refpp{fg1} where $F$ is perfect (see \refpp{fg1.0} for  details on the general case).

(i) $Aut(a_2(F)/Z(a_2(F)) \cong G_2(F)$, for  $ char(F) = 3$.

(ii) $Aut(g_2(F)) \cong  B_3(F)$, for $char(F) = 2$. Note that $Z(g_2(F)) = 0$ in this case.

(iii) $Aut(d_4(F)/Z(d_4(F)) \cong F_4(F)$, for $char(F) = 2$.

(iv) $Aut(d_n(F)/Z(d_n(F)) \cong B_n (F)$, for $char(F) = 2$, $n = 3$ or $n > 4$.
\end{enumerate}
 \end{prop}

\begin{rem}\labtt{solvable}   Line 2 of the table \refpp{XSg} is not  relevant here since we are assuming that $\frak g$ is based on an indecomposable root system.)  If $char(F)=2$ and $X=A_1$, then $\frak g$ is solvable and the reduced algebra is 2-dimensional abelian.   It automorphism group is $GL(2,F)$.
\end{rem}

\pf See \cite{steinberg} and \cite[page 210ff]{carter}. For the exceptional
cases, note that the automorphism group of the quotient algebra has no outer diagonal automorphisms, so contains all diagonal
automorphisms. \eop

\begin{thm} \labtt{autclassicalva}
Let $F$ be a field.
Let $V$ be either a VA of classical type over $F$ of non-exceptional type or
a reduced VA of classical type over $F$
in the exceptional case.

Then the restriction of
$Aut(V)$ to $V_1$ gives an isomorphism onto $Aut(V_1)$, the automorphism group of the reduced classical Lie algebra at degree 1 \refpp{autclassicalliealg}.
\end{thm}

\pf Obviously, the restriction is a monomorphism.   It suffices to show that the types of automorphisms indicated in
\refpp{autclassicalliealg} are represented.     By \refpp{chevgrouponva1}, \refpp{chevgrouponva2} and \refpp{largeaut}, we have
the Chevalley group $G(Y,F)$ for an appropriate type, $Y$,   and graph automorphisms in case $Y$ has type ADE; if so, the  standard integral form in $V_X$
admits graph automorphisms.
For arbitrary $Y$, it suffices to show that all diagonal automorphisms are represented on the VA.

Suppose that $Y$ has type ADE.   Let $L$ be the root lattice.
Think of the basis of $V_L$ over $F$ consisting of all $p\otimes e^{\a}$, where $p$ is a monomial.
A diagonal automorphism would be represented by some homomorphism $\l : Q \rightarrow F^{\times}$, where $Q$ is the root lattice.
The linear transformation defined by
$p\otimes e^{\a} \mapsto   \l (\a) \cdot p\otimes e^{\a}$
gives an automorphism of $V_L$ over $F$
corresponding to $\l$.

Suppose that $Y$ has type BCFG and we are not in an exceptional case.   Then $Y$ represents fixed points for action on a root system $X$ of type ADE by a graph automorphism $\g$.  Let $L$ be the root lattice for $X$.
A diagonal automorphism on $V(Y)$  would be represented by some homomorphism on the fixed point sublattice $\l : L^{\g} \rightarrow F^{\times}$.
Since $L^{\g}$ is a direct summand of $L$, the restriction map
$\dual L \rightarrow \dual {(L^{\g})}$ is onto.   Therefore, $\l $ extends to a homomorphism $\mu : L  \rightarrow F^{\times}$.   Then $\mu$ gives a diagonal outer automorphism on $V(X)$ whose restriction to $V(Y)$ is given by $\l$.

Suppose we are in an exceptional case.  The group of \refpp{autclassicalliealg}(ii), contains all diagonal automorphisms, so we are done.
\eop

\section{Modular moonshine and an embedding of the sporadic group $F_3$ into $E_8(3)$}

The simple group $F_3$ of order $2^{15}3^{10}5^3 7^2 13\cdot 19\cdot 31$ embeds in $E_8(3)$, a result of Peter Smith and John Thompson proved by computer work on a 248-dimensional lattice \cite{th}.   Richard Borcherds and Alex Ryba \cite{BR} suggested a vertex algebra style proof, but their program was not fully verified.

In this section, we show how a combination of the covering algebra (Section 4 and Appendix A)  and VOA integral form viewpoints allow us to complete the Borcherds-Ryba program.   It is possible that some other embeddings of finite groups into groups of Lie type could be proved with similar VOA methods.

\subsection{Modular moonshine of Borcherds and Ryba}\label{sec:6.1}

We mention one result of \cite{BR}.

\begin{thm}[{\cite[Corollary 4.8]{BR}}]\label{ch3c}
Let $g$ be a $3C$ element of the Monster group. Then there is a vertex algebra
$^gV=\oplus_{n\in \ZZ}  V_n$ defined over $\FF_3$  such that

 (1) $^gV$ is acted on by the group $C_\MM(g)/\la g\ra\cong F_3$; and

 (2) the Brauer trace  $\sum_{n\in \ZZ} \tr(h|^gV_n)q^{n-1} $ of any
3-regular element $h$ of $C_\MM(g)/\la g\ra$  on $^gV$ is a
Hauptmodul, and is equal to
the Hauptmodul of the element $gh$ of the monster. In particular, the
character of $^gV$ is equal to
\begin{equation*}
ch_{^gV}(q)= \sum_{n\in \ZZ} \dim(^gV_n)q^{n-1} = ch_{V_{E_8}}(q^3) =
q^{-1}+ 248q^2 + 4124q^5 +\cdots.
\end{equation*}
\end{thm}

Borcherds and Ryba also conjectured that $^gV$ (up to a certain rescaling of weights) is isomorphic to the vertex algebra
$IV_{E_8}/3IV_{E_8}$ over $\FF_3$.  (As discussed earlier in this article, $IV_{E_8}$ means the standard integral form for $V_{E_8}$ \refpp{standardif}.)  In this section, we shall give a proof of their conjecture.

\subsection{Lattice VOA $V_{E_8}$, Moonshine VOA and  $V_{EE_8}^+$}

In \cite{Mi}, Miyamoto gave a construction of the Moonshine VOA $V^\natural$ using 3 copies of the lattice type VOA
$V_{EE_8}^+$. An element  $g\in Aut(V^\natural)$ of the conjugacy class $3C$ is also constructed explicitly.  An alternative
treatment using the language of quadratic spaces was given by Shimakura in \cite{Shi}.  We shall  recall their construction and the definition of
a $3C$ element $g$ in this subsection.

\medskip

\begin{nota}\label{not6.1}
Let  $EE_8$  be the $\sqrt{2}$ times of the standard root lattice of type $E_8$. We will use $V_{EE_8}^+$ to denote the fixed point
subVOA of the lattice $V_{EE_8}$ by an involution $\theta\in Aut(V_{EE_8})$, which is a lift of the $-1$ isometry of $EE_8$.
\end{nota}

The representation theory for $V_{EE_8}^+$ has been studied by many authors \cite{ADL,Mi,Shi}. In particular, the set of all
irreducible modules and the fusion rules are determined.

\begin{prop}[\cite{ADL,Shi}]
The VOA $V_{EE_8}^+$ has exactly $2^{10}$ inequivalent irreducible modules and all of them are simple current modules.
Moreover,  the set of all inequivalent irreducible modules $R(V_{EE_8}^+)$ forms a 10-dimensional non-degenerate quadratic
space of  plus type over $\FF_2$ with the sum defined by fusion rules and a quadratic form given by the conformal weights (in
$\frac{1}2 \ZZ/\ZZ$).
\end{prop}

\begin{nota}
Let $M$ be an irreducible module of a VOA $V$. We use $[M]$ to denote the isomorphism class of $M$.
 \end{nota}

\begin{thm}\label{E8sum}
Let $\mathcal{S}$ be a maximal totally isotropic subspace of $R(V_{EE_8}^+)$. Then
\[
V(\mathcal{S}) =\bigoplus_{[M] \in \mathcal{S}} M
\]
has a unique structure of a  VOA
which extends the structure of the $M$ as  modules for  the summand $V_{EE_8}^+$.   Furthermore, this VOA is holomorphic and is of rank $8$ and so is isomorphic to $V_{E_8}$, the lattice VOA associated to the root lattice of type $E_8$.
\end{thm}

Let $\Phi$ and $\Psi$ be maximal totally singular subspaces of $R(V_{EE_8}^+)$ such that $\Phi\cap \Psi = 0$. For any $\{i, j\}
\subset \{1,2,3\}$ and set
\[
\Phi_{(i,j)} = \{(a_1, a_2,a_3) \in R(V_{EE_8}^+)^3\mid a_i = a_j \in \Phi,\text{ and } a_k = 0 \text{ if } k\notin \{i,j\}\}
\]
and
\[
\Psi_{(1,2,3)} = \{(b, b , b)\in  R(V_{EE_8}^+)^3\mid b\in \Psi\}.
\]

\begin{thm}[\cite{Shi}]\label{trans}
Let $ \mathcal{S}(\Phi,\Psi) = Span_{\FF_2} \{\Phi_{(1,2)}, \Phi_{(1,3)}, \Psi_{(1,2,3)} \}.$ Then  $\mathcal{S}(\Phi,\Psi)$ is a
maximal totally isotropic subspace of $R(V_{EE_8}^+)^3$ and
\[
V^\natural \cong V(\mathcal{S}(\Phi,\Psi)) =\bigoplus_{[A\otimes B\otimes C] \in \mathcal{S}(\Phi,\Psi)} A\otimes B\otimes C.
\]
as modules for $(V_{EE_8}^+)^3$.
Moreover, the automorphism group of $V^\natural$ is transitive on the set of all full subVOAs isomorphic to $(V_{EE_8}^+)^3$.
\end{thm}

\begin{nota}\labtt{defg}
Let $g$ be an automorphism of $V_{EE_8}^+ \otimes  V_{EE_8}^+ \otimes V_{EE_8}^+$ given by the cyclic permutation of the $3$
tensor factors.
 In \cite{Mi,Shi}, it was shown that $g$ can be lifted to an automorphism of $V^\natural$.   The lift is in the conjugacy class $3C$ of   $Aut(V^\natural)$.
\end{nota}
\medskip

\subsection{Positive definite real form of lattice VOAs}

The Moonshine VOA $V^\natural$ constructed in \cite{Mi} has a real form whose invariant bilinear form is positive definite. Next, we shall review the construction of a real form $\tilde{V}_{L,\RR}$ of lattice VOA $V_L$ such that the invariant form on $\tilde{V}_{L,\RR}$ is positive definite \cite{Mi}.

Recall that the lattice VOA constructed as in \cite{FLM} can be defined over $\RR$.
Let $V_{L,\RR}= S(\hat{\mathfrak{h}}^-_\RR)\otimes \RR\{L\}$ be the real lattice VOA associated to an even positive definite lattice, where $\mathfrak{h}_{\RR} =\RR\otimes_\ZZ L$, $ \hat{\mathfrak{h}}_{\RR}^-= \oplus_{n\in \ZZ^+} \mathfrak{h}\otimes \RR t^{-n}$. As usual, we use $x(-n)$ to denote $x\otimes t^{-n}$ for $x\in \mathfrak{h}$ and $n\in \ZZ^+$. As in \cite[Chapter 7]{FLM}, we fix a bilinear cocyle $\varepsilon(\ , \ )$ such that  $\varepsilon(\a,\a) = (-1)^{\la \a, \a\ra/2} $ and $\varepsilon( \a, \b)\epsilon (\b, \a) = (-1)^{\la \a, \b\ra}$.

\begin{nota}\labtt{liftof-1}
Let $\theta: V_{L,\RR} \to V_{L, \RR}$ be defined by
\[
\theta( x(-n_1)\cdots x(-n_k)\otimes e^\a) = (-1)^{k} x(-n_1)\cdots x(-n_k)\otimes e^{-\a}.
\]
Then $\theta$ is an automorphism of $V_{L,\RR}$, which is a lift of the $(-1)$-isometry of $L$ \cite{FLM}. We shall denote the $(\pm 1)$-eigenspaces of $\theta$ on $V_{L,\RR}$ by $V_{L,\RR}^\pm$.
\end{nota}

The following result is proved in \cite{Mi} (see also \cite{FLM}).
\begin{prop}[{\cite[Proposition 2.7]{Mi}}]\labtt{VLpd}
Let $L$ be an even positive definite lattice. Then
the real subspace
$\tilde{V}_{L,\RR}=V_{L,\RR}^+\oplus \sqrt{-1}V_{L,\RR}^-$
is a real form of $V_L$.
Furthermore,
the invariant form on $\tilde{V}_{L,\RR}$ is positive definite.
\end{prop}

\begin{rem}\label{twoauto}
Let $\varpi: V_L\to V_L$ be the anti-linear  involution induced by complex conjugation.
Then the real form $\tilde{V}_{L, \RR}$ defined in \refpp{VLpd} is in fact equal to the fixed point
subspace of $\theta\varpi$ in $V_L$.
\end{rem}

\begin{rem}\label{E8basis}
If $L=E_8$, then the Lie algebra of $(\tilde{V}_{E_8,\RR})_1$ is the compact real form of the  complex Lie algebra of type $E_8$.
For each positive root, let
\[
H_\alpha =\sqrt{-1} \alpha(-1),  \quad X_\alpha^+=e^\alpha + e^{-\alpha},  \quad
X_\alpha^-=\sqrt{-1}(e^\alpha -e^{-\alpha}).
\]
Then the Lie brackets are given as
\[
[H_\alpha, X^\pm_\beta] = (\alpha, \beta) X^\mp_\beta,
\]

\[
[X^\pm_\alpha, X^\pm_\beta]=
\begin{cases}
\pm \varepsilon(\alpha, \beta)  X^+_{\alpha+\beta}  &\text{ if } (\alpha, \beta)=-1,\\
0 & \text{otherwise},
\end{cases}
\]
and
\[
[X^+_\alpha, X^-_\beta]=
\begin{cases}
2H_\alpha  &\text{ if } (\alpha, \beta)=2,\\
\varepsilon(\alpha, \beta)  X^-_{\alpha+\beta}  &\text{ if } (\alpha, \beta)=-1,\\
0 & \text{otherwise}.
\end{cases}
\]
Let  $\{\a_1, \dots, \a_8\}$  be a set of simple roots of $E_8$. Then
\[
\{H_{\a_1},  \dots, H_{\a_8}\} \cup \{ X_\a^+, X_\a^-\mid \a \text{ positive roots}\},
\]
forms a basis of  $(\tilde{V}_{E_8,\RR})_1$.
\end{rem}

\subsection{Integral form and $\ZZ[\frac{1}2]$-form}\label{sec:6.3}

Let $IV_{E_8}$  and $IV_\Lambda$ be the standard integral forms of $V_{E_8}$ and $V_\Lambda$  \refpp{standardif}. In
\cite{BR}, it was shown that the integral form $IV_L$ contains the Virasoro element of $V_L$ if the lattice $L$ is unimodular. Since
$E_8$ and $\Lambda$ are unimodular lattices, the integral forms $IV_{E_8}$  and $IV_\Lambda$ also contain the Virasoro elements.

\medskip

Let $\theta$ be defined as in \refpp{liftof-1}. Then $\theta$ leave invariant the standard integral form $IV_{E_8}$ and $IV_{\Lambda}$.

\begin{nota}\labtt{12form}
Let $L$ be an even positive definite lattice and let $\mathcal{R}=\ZZ[1/2, \sqrt{-1}]$. Set  $\mathcal{R}V_L := IV_L\otimes_\ZZ
\mathcal{R}$ and define
\[
\tilde{V}_{L}[1/2] := (\mathcal{R}V_L)^{\la \theta\varpi\ra}= V_L[1/2]^+ \oplus  \sqrt{-1}V_L[1/2]^-,
\]
where $V_L[1/2] = IV_L\otimes_\ZZ \ZZ[1/2]$ and $V_L[1/2]^{\pm}$ is $\pm 1$-eigenspace of $\theta$ in $V_L[1/2]$.

Note that  $\tilde{V}_{L}[1/2]$ is a $\ZZ[1/2]$-form of $V_{L}$ and the invariant form on  $\tilde{V}_{L}[1/2]$ is positive definite.
\end{nota}

Now let $L=\Lambda$ be the Leech lattice and define
\[
V^0=V_{\Lambda,\RR}^+\cap \tilde{V}_\Lambda[1/2] = V_\Lambda[1/2]^+.
\]
Then $V^0$ is $\ZZ[1/2]$-form of $V_\Lambda^+$ and the invariant form on $V^0$ is positive definite. By using the triality automorphism defined in \cite{FLM}, one can also obtain a
$\ZZ[1/2]$-form of the Moonshine VOA $V^\natural = V_\Lambda^+\oplus V_\Lambda^{T,+}$ \cite[Theorem 3.2]{BR}.

\begin{prop}[\cite{BR}]\label{p6.7}
There is a $\ZZ[1/2]$-form $V^\natural [1/2]$ of the Moonshine VOA with the Virasoro element and a self-dual bilinear form.
\end{prop}

\begin{rem}
By the construction, the invariant form on $ V^\natural [1/2]$ is also positive definite.
\end{rem}

It is also well known \cite{CS,POE,LM} that the Leech lattice $\Lambda$ contains a sublattice isometric to $EE_8\perp EE_8\perp
EE_8$. Fix an embedding $\psi$ of
\[
EE_8\perp EE_8\perp EE_8 \hookrightarrow \Lambda.
\]
It induces an embedding of  $ (V_{EE_8}^+)^3$ into the Moonshine VOA $V^\natural$.  It is also shown in \cite{Shi} that the
subVOA $ (V_{EE_8}^+)^3$ defines an elementary abelian group $P$ of order $2^{15}$ such that
\[
(V^\natural)^P =  (V_{EE_8}^+)^3.
\]

Let $Irr(P)$ be the set of all irreducible characters of $P$ and set
\[
V^{\natural, \mu}= \{ v\in V^\natural \mid a\cdot v=\mu(a) v \text{ for any } a\in P\}, \quad \mu\in Irr(P).
\]
Define
\[
\mathcal{V}^{\natural, \mu}=V^{\natural, \mu}\cap V^{\natural}[1/2]\quad  \text{ for } \mu \in Irr(P).
\]
Then $\mathcal{V}^{\natural, \mu}$ is a $\ZZ[1/2]$-form of $V^{\natural, \mu}$ and $V^{\natural}[1/2] = \oplus_{\mu \in Irr(P)}
\mathcal{V}^{\natural, \mu}$.

\medskip

By \refpp{trans}, the automorphism group of $V^\natural$ is transitive on the set of all full subVOAs isomorphic to
$(V_{EE_8}^+)^3$.  Thus, using the above embedding $\psi$,
we have maximal totally singular subspaces $\Phi$ and $\Psi$ of $R(V_{EE_8}^+)$ such that
\[
V^\natural =   \bigoplus_{[A\otimes B\otimes C] \in \mathcal{S}(\Phi,\Psi)} A\otimes B\otimes C
\]
as modules for $(V_{EE_8}^+)^3$, where $ \mathcal{S}(\Phi,\Psi)$ is defined as in \refpp{trans}.

By \refpp{E8sum}, we also have
\[
V_{E_8} \cong \bigoplus_{[A]\in \Psi}  A
\]
as modules for $V_{EE_8}^+$.
Let $\tilde{V}_{E_8}[1/2]$ be the $\ZZ[1/2]$ form of $V_{E_8}$ as defined in \eqref{12form} and set $$A[1/2] =A\cap
\tilde{V}_{E_8}[1/2]\quad \text{ for }  [A]\in \Psi.$$

Now by \refpp{trans}, we have the following lemma.
\begin{lem}\label{embed}
Let  $\tilde{V}_{E_8}[1/2]$ and $A[1/2], [A]\in \Psi,$ be defined as above. Then there is  a monic linear transformation
\[
\eta: \tilde{V}_{E_8}[1/2] = \bigoplus_{[A]\in \Psi} A[1/2] \to    \bigoplus_{[A]\in \Psi} A[1/2]  \otimes A[1/2]\otimes A[1/2] \subset V^\natural[1/2]
\]
defined by
\[
\eta(x) =x\otimes x\otimes x  \quad \text{ for } x\in A[1/2], [A]\in \Psi.
\]
\end{lem}

The map $\eta$ will be used to create a covering transversal for the proof of \refpp{gv=ve83}.

\subsection{$\tilde{V}_{E_8}[1/2]$ modulo $3$}

Next we study the VA $\tilde{V}_{E_8}[1/2]/ 3\tilde{V}_{E_8}[1/2]$ over $\FF_3$.  First we recall a theorem of Lang \cite{Lang}
(see also \cite[E -Theorem 2.2]{Borel}).

\begin{thm}\label{langthm}
Let $G$ be a connected linear algebraic group and $\sigma$ an endomorphism of $G$ onto $G$ such that the fixed point
subgroup $G_\sigma$ is finite. Then the map $f : x\to  x\sigma(x)^{-1}$  of $G$ into $G$ is surjective.
\end{thm}

\begin{coro}\label{conjugate}  Let $q$ be a prime power.
If $G$ is an algebraic group over an algebraic extension of $\FF_q$, $\sigma$ is the Frobenius endomorphism based on $x
\mapsto x^q$, and if $\tau$ is an endomorphism of $G$ so that $\tau = r \sigma$, where $r$ is an inner automorphism, then
$\tau$ and $\sigma$ are conjugate.
\end{coro}
\pf
By \refpp{langthm}, there is $y$ in $G$ so that $ r=y \sigma (y)^{-1}$.  Then as functions from $G$ to $G$, we get
\[
y  \sigma y ^{-1} = y ( \sigma y ^{-1} \sigma^{-1}) \sigma =
y ( \sigma (y) ) ^{-1} \sigma =  r \sigma = \tau
\]
and the result follows.
\eop

\medskip

\begin{prop}\label{classicalE8}
Let $\tilde{V}_{E_8}[1/2]$ be the $\ZZ[1/2]$ form of $V_{E_8}$ as defined in \eqref{12form}. Then $\left(\tilde{V}_{E_8}[1/2]/
3\tilde{V}_{E_8}[1/2]\right)_1$ is isomorphic to a classical Lie algebra of type $E_8$ over $\FF_3$.
\end{prop}

\pf
Let $\mathcal{R}=\ZZ[1/2, \sqrt{-1}]$ be defined as in \refpp{12form}. Then
\[
\bar{\mathcal{R}}=\mathcal{R}/3\mathcal{R}\cong \FF_{9}.
\]
In this case, the anti-linear map $\varpi$ agrees with the Frobenius endomorphism $x\to x^3$ on $\bar{\mathcal{R}}$. Note that
$(a+b\sqrt{-1})^3 = a^3- b^3\sqrt{-1} = a-b\sqrt{-1}=\varpi(a+b\sqrt{-1})$.

Let $L= \left(\tilde{V}_{E_8}[1/2]/ 3\tilde{V}_{E_8}[1/2]\right)_1$. Then $\bar{\mathcal{R}} \otimes_{\FF_3} L$ is the classical
$E_8$ Lie algebra over $\FF_9$.

By \refpp{twoauto}, we have two automorphisms here. One is $\varpi$ which has fixed points equal to $M$,  the $\FF_3$ span of
the Chevalley basis in $\bar{\mathcal{R}} \otimes_{\FF_3} L$ and the other one is $\tau = \theta \varpi$ which has fixed points
$1\otimes L$.

By \refpp{conjugate},  these automorphisms are conjugate. Since the fixed points of one of them forms a classical $E_8$ algebra
over $\FF_3$, the same is true of the other. Therefore,  $L= \left(\tilde{V}_{E_8}[1/2]/ 3\tilde{V}_{E_8}[1/2]\right)_1$ is
isomorphic to  the classical $E_8$ algebra over $\FF_3$ and its automorphism group is $E_8(3)$.
\eop

\subsection{Definition of $^gV$}
In \cite{BR}, the vertex algebra $^gV$ is defined, as in Proposition \refpp{defgv} below.

\begin{prop}[{\cite[Corollary 4.8]{BR}}]\labtt{defgv}
Let $V^\natural [1/2]$ be the $\ZZ[1/2]$ form defined in \refpp{p6.7} (see also \cite{BR}).  Let $g$ be a $3C$-element of
the Monster. Denote
\[
{}^gV:= \frac{(V^\natural[1/2])^g}{ \nu(V^\natural[1/2])}.
\]
Then $^gV$ has the structure of a vertex algebra over  $\FF_3$.
\end{prop}

\begin{rem}\labtt{gvquotient} (i) The quotient ${}^gV= \frac{(V^\natural[1/2])^g}{ \nu(V^\natural[1/2])}$ is referred to as a 0-th Tate cohomology group in \cite{BR}.

(ii)
In \cite[Lemma 2.1 and Proposition 6.3]{BR},  it was shown that there is a non-singular invariant bilinear form on $^gV$  and
\[
^gV \cong \frac{V^{\la g \ra}  /3V^{\la g \ra}} {(V^{\la g \ra}  /3V^{\la g \ra}) \cap (V^{\la g \ra}  /3V^{\la g \ra})^\perp },
\]
where $V^{\la g\ra}$ is the fixed point subspace of $g$ in $V^\natural[1/2]$.
\end{rem}

\begin{lem}\label{Lie}
Let $V=\oplus_{n\geq 0} V_n$ be a VA over a field of characteristic $3$ such that
$\dim (V_0)=1$ and $V_n=0$ if $n\not\equiv 0 \mod 3$.  Then the weight 3 subspace $V_3$ has
a Lie algebra structure  given by the bracket
\[
[a, b] = a_2b  \quad \text{  for any  } a,b\in {V_3}
\]
and a symmetric invariant bilinear form
\[
(a|b )\cdot 1  =  a_5 b.
\]
\end{lem}

\pf
First we note that  $L(-1) a= a_{-2}\cdot 1 $ and hence $L(-1)$ is well defined. By skew-symmetry,
\[
[a,b]=a_2b = -b_2a +L(-1)  x,\quad \text{ for some } x\in {V_2}.
\]
Thus, we have $[a,b]=-[b,a]$ since $V_2=0$.

Moreover, by the commutator formula  (see Definition \ref{defva} (5) with $q=0$), we have

\[
[a_2, b_2] =\sum_{i=0}^2  \binom{2}{i}  (a_ib)_{2+2-i}.
\]
Since $V_4= V_5=0$, we have $a_1b=0$ and $a_0b=0$ and hence $[ a_2,
b_2]=(a_2b)_2$. Therefore, the  Jacobi Identity for
Lie algebra is satisfied.

Again by the skew symmetry, we have
\[
a_5b = (-1)^{1+5} b_5a =b_5a
\]
and hence $(a|b)=(b|a)$.

Finally, for any $a,b,c\in {V_3}$, we have
\[
a_5(b_2 c) = [a_5, b_2] c = \sum_{i=0}^5 \binom{5}{i} (a_i b)_{5+2-i} c.
\]
Note that $b_2(a_5c)=0$ because $b_2 1=0$ and $a_5c=(a|c) 1$.

Recall that $V_n=0$ if $n\not\equiv 0\mod 3$. Therefore, for any $a,b\in V$, $a_ib=0$ unless $i\equiv 2\mod 3$. Thus,  we have
\[
a_5(b_2 c) = 10 (a_2 b)_{5} c + (a_5b)_2 c  = (a_2 b)_5 c \mod 3.
\]
Hence, we have $(a|[b,c])= ([a,b]| c)$ and the form is invariant.
Note also that $(a_5b)_2 c= ((a|b)1)_2 c=0$.
\eop

\begin{lem}\label{affine}
Let $V=\oplus_{n\geq 0} V_n$ be a VA over a field of characteristic $3$ such that
$\dim V_0=1$ and $V_n=0$ if $n\not\equiv 0 \mod 3$.  Then for any $a, b\in V_3$, we have
\[
[a_{3m+2}, b_{3n+2}] = [a,b]_{3(m+n)+2} + m (a|b) \delta_{m+n,0} \quad \mod 3.
\]
\end{lem}

\pf
First we note that $a_i=0$ unless $i\equiv 2\mod 3$. By the commutator formula, we have
\[
[a_{3m+2}, b_{3n+2}] =\sum_{i=0}^{3m+2} \binom{3m+2}{i} (a_ib)_{3(m+n)+4-i}.
\]
Note that $a_ib=0$ unless $i=2$ or $5$. Moreover,
$\binom{3m+2}{2} \equiv 1 \mod 3$ and $\binom{3m+2}{5} \equiv m \mod 3$. Hence
we have
\[
[a_{3m+2}, b_{3n+2}] = [a,b]_{3(m+n)+2} + m (a|b) \delta_{m+n,0} \quad \mod 3
\]
as desired.
\eop

\begin{lem}
Let $^gV$ be defined as in \refpp{defgv}. Then the weight 3 subspace $^gV_3$ has
a Lie algebra structure  with a symmetric invariant bilinear form. Moreover, for any $a, b\in {^gV_3}$, we have
\[
[a_{3m+2}, b_{3n+2}] = [a,b]_{3(m+n)+2} + m (a|b) \delta_{m+n,0} \quad \mod 3.
\]
\end{lem}

\pf
The lemma follows from \refpp{Lie} and \refpp{affine} since ${^gV_n}=0$ if $n\not\equiv 0\mod 3$.
\eop

\subsection{$Im(\eta )$ provides a covering transversal}

\begin{lem}\label{eta}
Let $\eta: \tilde{V}_{E_8}[1/2] \to \bigoplus_{[A]\in \Psi} A[1/2]\otimes A[1/2]\otimes A[1/2]$ be the linear transformation defined as  in  \refpp{embed} and let $g$ be the $3C$ element defined as in  \refpp{defg} and let $\nu$ be the norm element for $\la g \ra$. Then
\[
(V^\natural[1/2])^g = \eta( \tilde{V}_{E_8}[1/2]) + \nu (V^\natural[1/2]).
\]
In particular, $\eta$ induces a linear surjection from
\[
 \tilde{V}_{E_8}[1/2] \to {^gV}= \frac{(V^\natural[1/2])^g}{ \nu(V^\natural[1/2])}.
\]
For simplicity, we also use $\eta$ to denote this surjection.
\end{lem}

\pf
By \refpp{trans} and the definition of $g$, we have
\[
(V^\natural[1/2])^g = \left(\bigoplus_{[A]\in \Psi} A[1/2]\otimes A[1/2]\otimes A[1/2]\right)^g + \nu(V^\natural[1/2]).
\]
Moreover,
\[
\begin{split}
&\  \left(\bigoplus_{[A]\in \Psi} A[1/2]\otimes A[1/2]\otimes A[1/2]\right)^g\\
  = & span_\ZZ \{ x\otimes x\otimes x\mid x\in A[1/2], A\in \Psi\}   + \nu(\bigoplus_{[A]\in \Psi} A[1/2]\otimes A[1/2]\otimes A[1/2]) \\
  = &  \eta( \tilde{V}_{E_8}[1/2]) + \nu(\bigoplus_{[A]\in \Psi} A[1/2]\otimes A[1/2]\otimes A[1/2]).
\end{split}
\]
Hence we have
\[
(V^\natural[1/2])^g = \eta( \tilde{V}_{E_8}[1/2]) + \nu(V^\natural[1/2])
\]
as desired.
\eop

\begin{de}
Let ${^gY}(\ , z) $ be the vertex operator of $^gV$. We define a new vertex operator
\[
{^gY}_{new}(u, z) ={^gY}(u, z^{1/3})\quad  \text{ for } u\in {^gV}.
\]
We also use $u_{\ell, new}$ to denote the coefficient of $z^{-\ell-1}$ in ${^gY}_{new}(u, z)$.
\end{de}

\begin{thm} \labtt{gv=ve83}
(i)
Let $\eta: \tilde{V}_{E_8}[1/2] \to  {^gV}$ be defined as in  \refpp{eta}, i.e,
\[
x\mapsto  x\otimes x\otimes x \quad  \text{ for } x\in A[1/2], [A]\in \Psi.
\]
Then $\eta$ induces  an isomorphism of VA between $\tilde{V}_{E_8}[1/2]/3\tilde{V}_{E_8}[1/2]$ and $(^gV, {^gY}_{new})$ over $\FF_3$. In
particular, $^gV_3$ is a Lie algebra of type $E_8$ over $\FF_3$.

\end{thm}

\pf   First, we recall  that $(a\otimes a\otimes a)_m (b\otimes b\otimes b)=0$ unless $m\equiv 2\mod 3$. Moreover, for any $\ell\in
\ZZ$, we have  $(a\otimes a\otimes a)_{\ell, new} = (a\otimes a\otimes a)_{3\ell+2}$.

For $[A], [B] \in \Psi $ and $\ell\in \ZZ$, and for any $a\in A[1/2], b\in B[1/2]$,
\[
(a\otimes a\otimes a)_{3\ell+2} (b\otimes b\otimes b) =
\sum_{i+j+k=3\ell} a_ib\otimes a_jb \otimes a_kb = a_\ell b\otimes a_\ell b \otimes a_\ell b.
\]
Therefore,  we have $\eta(a)_{\ell, new} \eta(b) = a_\ell b\otimes a_\ell b \otimes a_\ell b = \eta( a_\ell b)$. Thus, $\eta$ is a
homomorphism of VAs.  Since $\eta$ is clearly a surjection and by \eqref{ch3c},
\[
ch(\tilde{V}_{E_8}[1/2]/3\tilde{V}_{E_8}[1/2])= ch({^gV})= ch_{V_{E_8}}(q^3),
\]
$\eta$ is an isomorphism of VAs.

Now let  $a, b\in (\tilde{V}_{E_8}[1/2])_1$. Then
\[
(a\otimes a\otimes a)_2(b\otimes b\otimes b) = a_0b\otimes a_0b \otimes a_0b.
\]
Thus, $\eta$ also induces a Lie algebra isomorphism from $(\tilde{V}_{E_8}[1/2]/3\tilde{V}_{E_8}[1/2])_1$ to $({^gV})_3$.
Therefore, $({^gV})_3$ a Lie algebra of type $E_8$ over $\FF_3$ by \refpp{classicalE8}.
\eop

\begin{rem}\labtt{coveringtransversal2}  We give more detail about a covering transversal \refpp{coveringsubgroup} for a vertex algebra over $\FF_3$  in \refpp{gv=ve83}.
The vertex algebra is $V^{\natural} [\frac 12 ]^g / 3 V^{\natural} [\frac 12 ]^g$, the ideal is
$(\nu (V^{\natural} [\frac 12 ]^g)+3 V^{\natural} [\frac 12 ]^g )/3V^{\natural} [\frac 12 ]^g$
and the covering transversal is the subspace
$(Im(\eta )+3V^{\natural} [\frac 12 ]^g)/3V^{\natural} [\frac 12 ]^g$.   See \refpp{defgv}.
\end{rem}

\begin{coro} \labtt{f3ine83}   The sporadic simple group $F_3$, of order $2^{15}3^{10}5^3 7^2 13{\cdot }19{\cdot} 31$,  embeds in $E_8(3)$.
\end{coro}
\pf
The group $F_3$ may be defined as $C(g)/\la g\ra$,
where $g$ is a class $3C$-element of the Monster and $C(g)$ denotes the centralizer in the Monster of $g$.    We clearly have an
action of $C(g)/\la g\ra$ on $\vnat [1/2]$, whence an action on the space $^g V$ of \refpp{gvquotient}.    Now use the
isomorphism with $\tilde{V}_{E_8}[1/2]/3\tilde{V}_{E_8}[1/2]$,
whose automorphism group is $E_8(3)$.
\eop

\vspace{1cm}

\centerline{\bf Acknowledgments}

\smallskip

The first author thanks Academia Sinica for hospitality during visits in 2012 and 2013, and the US National Security Agency and the University of Michigan for financial support.  The second author thanks National Science Council (NSC 100-2628-M-001005-MY4) and National Center for Theoretical Sciences of Taiwan for financial support.

\newpage

\appendix
\centerline{\bf \Large Appendices}

\section{The covering algebra procedure for classical Lie algebras over  algebraically closed field \labtt{ca} }

See the article \cite{fg} for background.   A copy is posted at
\begin{verbatim}
http://www.math.lsa.umich.edu/~rlg/griesspublicationlist.html
\end{verbatim}

The main results of \cite{fg} are summarized at the end of this section,  \refpp{fg1}, \refpp{fg2}.   We first give a few definitions and elementary results on algebras and automorphisms.

\begin{de} \labtt{quasisimple}   If $X$ is a group or a Lie algebra, $X$ is quasisimple if $X$ modulo its center is simple.
\end{de}

We take two elementary results from \cite{fg}:

\begin{lem} \labtt{fglem1}
 Let $\s$ be an automorphism of an algebra over the field $K$. Suppose that $C$ is the fixed
point subalgebra and $\phi \in K[t]$. Then $N := Im(\phi (\s ))$ is stable under multiplication
by $C$. In particular, $C\cap N$ is an ideal of $C$.
 \end{lem}
\pf Straightforward. \eop

The following is an immediate consequence.

\begin{coro} \labtt{fglem2}  Let $A$ be an algebra over the field $K$, and let $\s$ be an automorphism
of $A$ with fixed point subalgebra $C$. Suppose that $\s$ has minimal polynomial $(t-1)^k$.   Define
$N:=Im((1-\s )^{k-1})$.
(i) If  $C\ne N$,  then $C$ is not simple.
(ii) If $N \not \le Z(C)$, then  $Z(C)$ is an ideal and $C \ne N + Z(C)$, then $C$ is not quasi-simple.
\end{coro}

\begin{de}\labtt{classicalliealg}
Let $X$ be a root system and $L$ the lattice spanned by a Chevalley basis in a Lie algebra of type $X$.   A Lie algebra over the commutative ring $R$ is said to have classical type $X$ if it is isomorphic to $R\otimes L$.
\end{de}

\begin{de} \labtt{exceptionalpair} Suppose that $X$ is an indecomposable root system of type ADE and $p$ is a prime number so that $X$ has a graph automorphism, say $\g$,  of order $p$.
Let $X'$ denote the set of roots fixed by $\g$.

Let $F$ be a field of characteristic $p$, $\frak g$ be a Lie algebra over $F$ of classical type $X$ and let $\frak g'$ be the subalgebra corresponding to $X'$.      Then the pair $(X', p)$ and subalgebra of type $X'$ corresponding to $\g$ are each called exceptional.
The central quotients of $\frak g'$ have exceptionally large automorphism groups, and all such cases are described in \refpp{fg1}.

Given an exceptional pair, $(X',p)$,
there is an essentially unique pair $(X,\g)$ which gives rise to $(X',p)$ as above.   We call such $(X,\g)$ the ancestor of $(X',p)$.  The quotient of the Lie algebra, resp. VA of type $X'$, by its norm ideal is called the reduced Lie algebra, resp. reduced VA of type $X'$.
\end{de}

For all exceptional cases, we use   polynomials in the graph automorphism to define a fixed point subalgebra and an ideal.
  The covering algebra \refpp{coveringsubgroup}  turns out to be $\frak g'$ in the notation of \refpp{exceptionalpair}.

For example, in  case (i) of the Theorem \refpp{fg1}, we start with Lie algebra of type $D_4$ in characteristic 3, $\g$ is a graph automorphism of order 3.  The fixed point subalgebra has type $g_2$ and in it, we have a 7-dimensional ideal and covering algebra of type $a_2$.

\begin{thm}\labtt{fg1}  If F is algebraically closed or perfect, then

(i) $Aut(a_2(F)/Z(a_2(F)) \cong G_2(F)$, for  $ char(F) = 3$.

(ii) $Aut(g_2(F)) \cong  B_3(F)$, for $char(F) = 2$. Note that $Z(g_2(F)) = 0$ in this case.

(iii) $Aut(d_4(F)/Z(d_4(F)) \cong F_4(F)$, for $char(F) = 2$.

(iv) $Aut(d_n(F)/Z(d_n(F)) \cong B_n (F)$, for $char(F) = 2$, $n = 3$ or $n > 4$.

If $F$ is arbitrary, then in each case, the left side contains a normal subgroup isomorphic to the right side, giving quotient an abelian group of exponent $char(F) > 0$.
\end{thm}

\begin{rem}\labtt{fg1.0}   The case where $F$ is algebraically closed was treated in \cite{fg}.   Suppose that $F$ is arbitrary, of characteristic  $p= $ 2 or 3.   Let $K$ be an algebraic closure and let $\varphi$ be the Frobenius endomorphism of the algebraic group $G$ indicated by one of the cases (i, ii, iii, iv).
Then by taking fixed points of $\varphi$, we deduce the conclusions (i, ii, iii, iv) over the prime field.   By extending the field to $F$, we get groups $G(F)$ over $F$ contained in  $Aut(\frak q (F) )$, where $\frak q (F)$ is  the quotient algebra indicated in the respective cases.
If $F$ is perfect, the containments $G(F) \le Aut(\frak q (F) )$ are equalities since $G(K) = Aut(\frak q (K) )$ and the groups $G(F)$ have no nontrivial diagonal outer automorphisms.
If $F$ is not perfect, the natural quotient $Aut(\frak q (F) )/G(F)$ is covered by the quotient of  $F^{\times}$ by its group of $p$-th powers.
\end{rem}

\begin{rem}\labtt{fg1.1}
For details of the proofs, see \cite{fg}.    In \refpp{fg1}, $\s$ is a standard graph automorphism of order $p=char(K)$ and
$N=Im(\sum_{i=0}^{p-1} \s^i )$ is contained in $C$, so corresponds to the Lie subalgebra of $C$ associated to the short roots.
\end{rem}

\begin{thm}\labtt{fg2}     The following Lie algebras are not quasi-simple:

Types $B_n(K), C_n(K), F_4(K)$ when $char(K) = 2$.

Type $G_2(K)$ when $char(K) = 3$.
\end{thm}

The latter result means that the proper nonzero ideal created by the procedure of \cite{fg} is not central.  Automorphism groups for quotients by various subspaces of the center may be found in \cite{fg}.

\end{document}